\documentclass{amsart}%
\usepackage{amsfonts}
\usepackage{graphicx}
\usepackage{amscd}%
\usepackage{amsmath}%
\usepackage{amssymb}
\newtheorem{theorem}{Theorem}
\theoremstyle{plain}

\newtheorem{proposition}{Proposition}

\numberwithin{equation}{section}
\numberwithin{theorem}{section}
\numberwithin{algorithm}{section}
\numberwithin{axiom}{section}
\numberwithin{case}{section}
\numberwithin{claim}{section}
\numberwithin{conclusion}{section}
\numberwithin{condition}{section}
\numberwithin{conjecture}{section}
\numberwithin{corollary}{section}
\numberwithin{criterion}{section}
\numberwithin{definition}{section}
\numberwithin{example}{section}
\numberwithin{exercise}{section}
\numberwithin{lemma}{section}
\numberwithin{notation}{section}
\numberwithin{problem}{section}
\numberwithin{proposition}{section}
\numberwithin{remark}{section}
\numberwithin{solution}{section}

\begin{document}
\title[On an integral system]{On the integral systems related to Hardy-Littlewood-Sobolev inequality}
\author{Fengbo Hang}
\address{Department of Mathematics\\
Princeton University\\
Fine Hall, Washington Road\\
Princeton, NJ 08544}
\email{fhang@math.princeton.edu}

\begin{abstract}
We prove all the maximizers of the sharp Hardy-Littlewood-Sobolev inequality
are smooth. More generally, we show all the nonnegative critical functions are
smooth, radial with respect to some points and strictly decreasing in the
radial direction. In particular, we resolve all the cases left open by
previous works of Chen, Li and Ou on the corresponding integral systems.

\end{abstract}
\maketitle

\section{Introduction\label{sec1}}

The classical Hardy-Littlewood-Sobolev inequality states that for $0<\alpha
<n$, $1<p_{0},q_{0}<\frac{n}{\alpha}$ such that $\frac{1}{p_{0}}+\frac
{1}{q_{0}}=1+\frac{\alpha}{n}$ (see \cite[theorem 1 on p119]{S})%
\[
\left\vert \int_{\mathbb{R}^{n}\times\mathbb{R}^{n}}\frac{f\left(  x\right)
g\left(  y\right)  }{\left\vert x-y\right\vert ^{n-\alpha}}dxdy\right\vert
\leq c\left(  n,p_{0},\alpha\right)  \left\vert f\right\vert _{L^{p_{0}%
}\left(  \mathbb{R}^{n}\right)  }\left\vert g\right\vert _{L^{q_{0}}\left(
\mathbb{R}^{n}\right)  }.
\]
In \cite{Li}, it was shown that the sharp constant%
\[
c\left(  n,p_{0},\alpha\right)  =\sup\left\{  \int_{\mathbb{R}^{n}%
\times\mathbb{R}^{n}}\frac{f\left(  x\right)  g\left(  y\right)  }{\left\vert
x-y\right\vert ^{n-\alpha}}dxdy:\left\vert f\right\vert _{L^{p_{0}}\left(
\mathbb{R}^{n}\right)  }=1,\left\vert g\right\vert _{L^{q_{0}}\left(
\mathbb{R}^{n}\right)  }=1\right\}
\]
is achieved by some functions $f$ and $g$. Moreover, after multiplying some
constants, any maximizer $f,g$ must be radial symmetric with respect to the
same point, strictly decreasing in the radial direction and satisfy the
integral system%
\[
f\left(  x\right)  ^{p_{0}-1}=\int_{\mathbb{R}^{n}}\frac{g\left(  y\right)
}{\left\vert x-y\right\vert ^{n-\alpha}}dy,\quad g\left(  x\right)  ^{q_{0}%
-1}=\int_{\mathbb{R}^{n}}\frac{f\left(  y\right)  }{\left\vert x-y\right\vert
^{n-\alpha}}dy.
\]
It was also shown that when $p_{0}=q_{0}$, we have%
\[
f\left(  x\right)  =g\left(  x\right)  =c\left(  n,p_{0}\right)  \left(
\frac{\lambda}{\left\vert x-x_{0}\right\vert ^{2}+\lambda^{2}}\right)
^{n/p_{0}}%
\]
for some $\lambda>0$ and $x_{0}\in\mathbb{R}^{n}$.

If we let $p=\frac{1}{p_{0}-1}$, $q=\frac{1}{q_{0}-1}$, $u=f^{p_{0}-1}$,
$v=g^{q_{0}-1}$, then the Euler-Lagrange equation becomes%
\begin{equation}
u\left(  x\right)  =\int_{\mathbb{R}^{n}}\frac{v\left(  y\right)  ^{q}%
}{\left\vert x-y\right\vert ^{n-\alpha}}dy,\quad v\left(  x\right)
=\int_{\mathbb{R}^{n}}\frac{u\left(  y\right)  ^{p}}{\left\vert x-y\right\vert
^{n-\alpha}}dy\label{eq1.1}%
\end{equation}
for nonnegative functions $u\in L^{p+1}\left(  \mathbb{R}^{n}\right)  $ and
$v\in L^{q+1}\left(  \mathbb{R}^{n}\right)  $ and $0<\alpha<n$, $\frac{\alpha
}{n-\alpha}<p,q<\infty$, $\frac{1}{p+1}+\frac{1}{q+1}+\frac{\alpha}{n}=1$.
When $p=q=\frac{n+\alpha}{n-\alpha}$, as observed in \cite{Li}, it follows
from the fact $\frac{1}{\left\vert x\right\vert ^{n-\alpha}}=c\left(
n,\alpha\right)  \frac{1}{\left\vert x\right\vert ^{n-\frac{\alpha}{2}}}%
\ast\frac{1}{\left\vert x\right\vert ^{n-\frac{\alpha}{2}}}$ that $u=v$, then
the system reduces to%
\begin{equation}
u\left(  x\right)  =\int_{\mathbb{R}^{n}}\frac{u\left(  y\right)
^{\frac{n+\alpha}{n-\alpha}}}{\left\vert x-y\right\vert ^{n-\alpha}%
}dy.\label{eq1.2}%
\end{equation}
In \cite{CLO2}, using an integral form of the method of moving planes
(\cite{GNN}), it was shown that any nonzero nonnegative regular solution $u$
of (\ref{eq1.2}) must be of the form%
\[
u\left(  x\right)  =c\left(  n,\alpha\right)  \left(  \frac{\lambda}%
{\lambda^{2}+\left\vert x-x_{0}\right\vert ^{2}}\right)  ^{\frac{n-\alpha}{2}}%
\]
for some $\lambda>0$ and $x_{0}\in\mathbb{R}^{n}$. This solves an open problem
proposed in \cite{Li} (see a somewhat different argument in \cite{L} and the
clarifications in \cite[Remark 1.3 on p332]{CLO2}). In \cite{CL,CLO1}, such
kinds of analysis were extended to the system (\ref{eq1.1}) under the
additional constraints $p\geq1$ and $q\geq1$. However the analysis does not
give the regularity of maximizer for all the Hardy-Littlewood-Sobolev
inequalities. On the other hand, it does not seem that we will have nonsmooth
maximizers for the Hardy-Littlewood-Sobolev inequality in any case. The main
aim of this article is to prove the regularity and radial symmetry of
nonnegative solutions of the system (\ref{eq1.1}) in its full range. Another
motivation comes from the study of regularity issues for a similar integral
system in \cite{HWY}.

\begin{theorem}
\label{thmmain}Assume $0<\alpha<n$, $\frac{\alpha}{n-\alpha}<p,q<\infty$,
$\frac{1}{p+1}+\frac{1}{q+1}+\frac{\alpha}{n}=1$, $u\in L^{p+1}\left(
\mathbb{R}^{n}\right)  $ is nonnegative and does not vanish identically. If%
\[
v\left(  x\right)  =\int_{\mathbb{R}^{n}}\frac{u\left(  y\right)  ^{p}%
}{\left\vert x-y\right\vert ^{n-\alpha}}dy,\quad u\left(  x\right)
=\int_{\mathbb{R}^{n}}\frac{v\left(  y\right)  ^{q}}{\left\vert x-y\right\vert
^{n-\alpha}}dy.
\]
Then $u\in C^{\infty}\left(  \mathbb{R}^{n}\right)  $, $v\in C^{\infty}\left(
\mathbb{R}^{n}\right)  $. Moreover, there exists a point $x_{0}\in
\mathbb{R}^{n}$ such that both $u$ and $v$ are radial symmetric with respect
to $x_{0}$ and strictly decreasing along radial direction.
\end{theorem}

Indeed, the regularity is still true under the relatively weaker assumption
$u\in L_{loc}^{p+1}\left(  \mathbb{R}^{n}\right)  $ (see Proposition
\ref{propregularity}). The method in \cite{CL,CLO1}, which is basically linear
in nature, does not seem to work for the case when one of the two indices $p$
and $q$ is strictly less than $1$. We will develop some nonlinear approaches
which work for all $p$ and $q$ at once. In \cite{HWY}, we will apply this
technique to derive the regularity for another integral system. In Section
\ref{sec2} below, we will prove a local regularity result which has the
regularity part in Theorem \ref{thmmain} as a corollary. In Section
\ref{sec3}, we will prove all the solutions are radial.

\textbf{Acknowledgment}: The research of the author is supported by National
Science Foundation Grant DMS-0501050 and a Sloan Research Fellowship. Part of
the work was done while I was visiting MSRI, I would like to thank the
institute for hospitality. I also thank the anonymous referee for his/her
suggestions which improve the presentation of the article.

\section{Regularity issue\label{sec2}}

In this section, we will show any solution $u,v$ to the system (\ref{eq1.1})
must be smooth if we assume $u\in L_{loc}^{p+1}\left(  \mathbb{R}^{n}\right)
$. Such a local integrability condition is necessary for the smoothness
because as observed in \cite{Li}, system (\ref{eq1.1}) has singular solutions
as%
\[
u\left(  x\right)  =c\left(  n,\alpha,p\right)  \left\vert x\right\vert
^{-\frac{n}{p+1}},\quad v\left(  x\right)  =c\left(  n,\alpha,p\right)
\left\vert x\right\vert ^{-\frac{n}{q+1}}.
\]
This follows from a simple change of variable in the integrals. To achieve the
regularity, we start with a local result which has some similarity to
\cite[theorem 2]{CJLL} and \cite[theorem 1.3]{L}.

\begin{proposition}
\label{proplocal}Given $0<\alpha,\beta<n$, $1<a,b\leq\infty$, $1\leq r<\infty$
such that%
\[
\frac{1}{ra}+\frac{1}{b}=\frac{\alpha}{rn}+\frac{\beta}{n}.
\]
Assume%
\[
\frac{n}{n-\beta}<p<q<\infty,
\]%
\[
\frac{\alpha}{n}<\frac{r}{q}+\frac{1}{a}<\frac{r}{p}+\frac{1}{a}<1,
\]
$u,f\in L^{p}\left(  B_{R}\right)  $, $U\in L^{a}\left(  B_{R}\right)  $,
$V\in L^{b}\left(  B_{R}\right)  $ are all nonnegative functions with $\left.
f\right\vert _{B_{R/2}}\in L^{q}\left(  B_{R/2}\right)  $,%
\[
\left\vert U\right\vert _{L^{a}\left(  B_{R}\right)  }^{1/r}\left\vert
V\right\vert _{L^{b}\left(  B_{R}\right)  }\leq\varepsilon\left(
n,p,q,r,\alpha,\beta,a,b\right)  \text{ small}%
\]
and%
\[
u\left(  x\right)  \leq\int_{B_{R}}\frac{V\left(  y\right)  }{\left\vert
x-y\right\vert ^{n-\beta}}\left[  \int_{B_{R}}\frac{U\left(  z\right)
u\left(  z\right)  ^{r}}{\left\vert y-z\right\vert ^{n-\alpha}}dz\right]
^{1/r}dy+f\left(  x\right)
\]
for $x\in B_{R}$, then $u\in L^{q}\left(  B_{R/4}\right)  $, moreover%
\[
\left\vert u\right\vert _{L^{q}\left(  B_{R/4}\right)  }\leq c\left(
n,p,q,r,\alpha,\beta,a,b\right)  \left(  R^{\frac{n}{q}-\frac{n}{p}}\left\vert
u\right\vert _{L^{p}\left(  B_{R}\right)  }+\left\vert f\right\vert
_{L^{q}\left(  B_{R/2}\right)  }\right)  .
\]

\end{proposition}

\begin{proof}
By scaling, we may assume $R=1$. First assume we have $u,f\in L^{q}\left(
B_{1}\right)  $. Denote%
\[
v\left(  x\right)  =\int_{B_{1}}\frac{U\left(  y\right)  u\left(  y\right)
^{r}}{\left\vert x-y\right\vert ^{n-\alpha}}dy\text{ for }x\in B_{1}.
\]
Let $p_{1}$ and $q_{1}$ be the numbers defined by%
\[
\frac{1}{p_{1}}=\frac{r}{p}+\frac{1}{a}-\frac{\alpha}{n},\quad\frac{1}{q_{1}%
}=\frac{r}{q}+\frac{1}{a}-\frac{\alpha}{n},
\]
then it follows from Hardy-Littlewood-Sobolev inequality that%
\begin{align*}
\left\vert v\right\vert _{L^{p_{1}}\left(  B_{1}\right)  }  & \leq c\left(
n,p,r,\alpha,a\right)  \left\vert U\right\vert _{L^{a}\left(  B_{1}\right)
}\left\vert u\right\vert _{L^{p}\left(  B_{1}\right)  }^{r},\\
\left\vert v\right\vert _{L^{q_{1}}\left(  B_{1}\right)  }  & \leq c\left(
n,q,r,\alpha,a\right)  \left\vert U\right\vert _{L^{a}\left(  B_{1}\right)
}\left\vert u\right\vert _{L^{q}\left(  B_{1}\right)  }^{r}.
\end{align*}
Given $0<s<t\leq1/2$. For $x\in B_{s}$, we have%
\begin{align*}
u\left(  x\right)   & \leq\int_{B_{\frac{s+t}{2}}}\frac{V\left(  y\right)
v\left(  y\right)  ^{1/r}}{\left\vert x-y\right\vert ^{n-\beta}}dy+\int
_{B_{1}\backslash B_{\frac{s+t}{2}}}\frac{V\left(  y\right)  v\left(
y\right)  ^{1/r}}{\left\vert x-y\right\vert ^{n-\beta}}dy+f\left(  x\right) \\
& \leq\int_{B_{\frac{s+t}{2}}}\frac{V\left(  y\right)  v\left(  y\right)
^{1/r}}{\left\vert x-y\right\vert ^{n-\beta}}dy+\frac{c\left(  n,\beta\right)
}{\left(  t-s\right)  ^{n-\beta}}\int_{B_{1}\backslash B_{\frac{s+t}{2}}%
}V\left(  y\right)  v\left(  y\right)  ^{1/r}dy+f\left(  x\right) \\
& \leq\int_{B_{\frac{s+t}{2}}}\frac{V\left(  y\right)  v\left(  y\right)
^{1/r}}{\left\vert x-y\right\vert ^{n-\beta}}dy+\frac{c\left(  n,p,r,\alpha
,\beta,a,b\right)  \left\vert u\right\vert _{L^{p}\left(  B_{1}\right)  }%
}{\left(  t-s\right)  ^{n-\beta}}+f\left(  x\right)  .
\end{align*}
Hence we have%
\[
\left\vert u\right\vert _{L^{q}\left(  B_{s}\right)  }\leq c\left(
n,q,r,\beta,b\right)  \left\vert V\right\vert _{L^{b}\left(  B_{1}\right)
}\left\vert v\right\vert _{L^{q_{1}}\left(  B_{\frac{s+t}{2}}\right)  }%
^{1/r}+\frac{c\left(  n,p,q,r,\alpha,\beta,a,b\right)  \left\vert u\right\vert
_{L^{p}\left(  B_{1}\right)  }}{\left(  t-s\right)  ^{n-\beta}}+\left\vert
f\right\vert _{L^{q}\left(  B_{1/2}\right)  }.
\]
On the other hand, for $x\in B_{\frac{s+t}{2}}$, we have%
\begin{align*}
v\left(  x\right)   & =\int_{B_{t}}\frac{U\left(  y\right)  u\left(  y\right)
^{r}}{\left\vert x-y\right\vert ^{n-\alpha}}dy+\int_{B_{1}\backslash B_{t}%
}\frac{U\left(  y\right)  u\left(  y\right)  ^{r}}{\left\vert x-y\right\vert
^{n-\alpha}}dy\\
& \leq\int_{B_{t}}\frac{U\left(  y\right)  u\left(  y\right)  ^{r}}{\left\vert
x-y\right\vert ^{n-\alpha}}dy+\frac{c\left(  n,\alpha\right)  }{\left(
s-t\right)  ^{n-\alpha}}\int_{B_{1}\backslash B_{t}}U\left(  y\right)
u\left(  y\right)  ^{r}dy\\
& \leq\int_{B_{t}}\frac{U\left(  y\right)  u\left(  y\right)  ^{r}}{\left\vert
x-y\right\vert ^{n-\alpha}}dy+\frac{c\left(  n,p,r,\alpha,a\right)  \left\vert
U\right\vert _{L^{a}\left(  B_{1}\right)  }\left\vert u\right\vert
_{L^{p}\left(  B_{1}\right)  }^{r}}{\left(  s-t\right)  ^{n-\alpha}}.
\end{align*}
This implies%
\[
\left\vert v\right\vert _{L^{q_{1}}\left(  B_{\frac{s+t}{2}}\right)  }\leq
c\left(  n,q,r,\alpha,a\right)  \left\vert U\right\vert _{L^{a}\left(
B_{1}\right)  }\left\vert u\right\vert _{L^{q}\left(  B_{t}\right)  }%
^{r}+\frac{c\left(  n,p,q,r,\alpha,a\right)  \left\vert U\right\vert
_{L^{a}\left(  B_{1}\right)  }\left\vert u\right\vert _{L^{p}\left(
B_{1}\right)  }^{r}}{\left(  s-t\right)  ^{n-\alpha}}.
\]
Combine the two inequalities together, we see%
\begin{align*}
\left\vert u\right\vert _{L^{q}\left(  B_{s}\right)  }  & \leq c\left(
n,q,r,\alpha,\beta,a,b\right)  \left\vert U\right\vert _{L^{a}\left(
B_{1}\right)  }^{1/r}\left\vert V\right\vert _{L^{b}\left(  B_{1}\right)
}\left\vert u\right\vert _{L^{q}\left(  B_{t}\right)  }\\
& +\frac{c\left(  n,p,q,r,\alpha,\beta,a,b\right)  }{\left(  s-t\right)
^{\max\left\{  \left(  n-\alpha\right)  /r,n-\beta\right\}  }}\left\vert
u\right\vert _{L^{p}\left(  B_{1}\right)  }+\left\vert f\right\vert
_{L^{q}\left(  B_{1/2}\right)  }\\
& \leq\frac{1}{2}\left\vert u\right\vert _{L^{q}\left(  B_{t}\right)  }%
+\frac{c\left(  n,p,q,r,\alpha,\beta,a,b\right)  }{\left(  s-t\right)
^{\max\left\{  \left(  n-\alpha\right)  /r,n-\beta\right\}  }}\left\vert
u\right\vert _{L^{p}\left(  B_{1}\right)  }+\left\vert f\right\vert
_{L^{q}\left(  B_{1/2}\right)  }%
\end{align*}
if $\varepsilon$ is small enough. It follows from usual iteration procedure
(\cite[lemma 4.3 on p.75]{HL}) that%
\[
\left\vert u\right\vert _{L^{q}\left(  B_{1/4}\right)  }\leq c\left(
n,p,q,r,\alpha,\beta,a,b\right)  \left(  \left\vert u\right\vert
_{L^{p}\left(  B_{1}\right)  }+\left\vert f\right\vert _{L^{q}\left(
B_{1/2}\right)  }\right)  .
\]
To prove the full proposition, we note that for some function $0\leq
\eta\left(  x\right)  \leq1$,%
\[
u\left(  x\right)  =\eta\left(  x\right)  \int_{B_{R}}\frac{V\left(  y\right)
}{\left\vert x-y\right\vert ^{n-\beta}}\left[  \int_{B_{R}}\frac{U\left(
z\right)  u\left(  z\right)  ^{r}}{\left\vert y-z\right\vert ^{n-\alpha}%
}dz\right]  ^{1/r}dy+\eta\left(  x\right)  f\left(  x\right)  .
\]
We may define a map $T$ by%
\[
T\left(  \varphi\right)  \left(  x\right)  =\eta\left(  x\right)  \int_{B_{R}%
}\frac{V\left(  y\right)  }{\left\vert x-y\right\vert ^{n-\beta}}\left[
\int_{B_{R}}\frac{U\left(  z\right)  \left\vert \varphi\left(  z\right)
\right\vert ^{r}}{\left\vert y-z\right\vert ^{n-\alpha}}dz\right]  ^{1/r}dy.
\]
Note that we have%
\[
\left\vert T\left(  \varphi\right)  \right\vert _{L^{p}\left(  B_{1}\right)
}\leq c\left(  n,p,r,\alpha,\beta,a,b\right)  \left\vert U\right\vert
_{L^{a}\left(  B_{1}\right)  }^{1/r}\left\vert V\right\vert _{L^{b}\left(
B_{1}\right)  }\left\vert \varphi\right\vert _{L^{p}\left(  B_{1}\right)
}\leq\frac{1}{2}\left\vert \varphi\right\vert _{L^{p}\left(  B_{1}\right)  }%
\]
and%
\[
\left\vert T\left(  \varphi\right)  \right\vert _{L^{q}\left(  B_{1}\right)
}\leq c\left(  n,q,r,\alpha,\beta,a,b\right)  \left\vert U\right\vert
_{L^{a}\left(  B_{1}\right)  }^{1/r}\left\vert V\right\vert _{L^{b}\left(
B_{1}\right)  }\left\vert \varphi\right\vert _{L^{q}\left(  B_{1}\right)
}\leq\frac{1}{2}\left\vert \varphi\right\vert _{L^{q}\left(  B_{1}\right)  }%
\]
if $\varepsilon$ is small enough. Moreover, for any $\varphi,\psi\in
L^{p}\left(  B_{1}\right)  $, it follows from Minkowski's inequality that%
\[
\left\vert T\left(  \varphi\right)  \left(  x\right)  -T\left(  \psi\right)
\left(  x\right)  \right\vert \leq T\left(  \left\vert \varphi-\psi\right\vert
\right)  \left(  x\right)  \text{ for }x\in B_{1},
\]
hence%
\[
\left\vert T\left(  \varphi\right)  -T\left(  \psi\right)  \right\vert
_{L^{p}\left(  B_{1}\right)  }\leq\left\vert T\left(  \left\vert \varphi
-\psi\right\vert \right)  \right\vert _{L^{p}\left(  B_{1}\right)  }\leq
\frac{1}{2}\left\vert \varphi-\psi\right\vert _{L^{p}\left(  B_{1}\right)  }.
\]
Similarly, we have for any $\varphi,\psi\in L^{q}\left(  B_{1}\right)  $,%
\[
\left\vert T\left(  \varphi\right)  -T\left(  \psi\right)  \right\vert
_{L^{q}\left(  B_{1}\right)  }\leq\frac{1}{2}\left\vert \varphi-\psi
\right\vert _{L^{q}\left(  B_{1}\right)  }.
\]
For $k\in\mathbb{N}$, let $f_{k}\left(  x\right)  =\min\left\{  f\left(
x\right)  ,k\right\}  $, then it follows from contraction mapping theorem that
we may find a unique $u_{k}\in L^{q}\left(  B_{1}\right)  $ such that%
\begin{align*}
u_{k}\left(  x\right)   & =T\left(  u_{k}\right)  \left(  x\right)
+\eta\left(  x\right)  f_{k}\left(  x\right) \\
& =\eta\left(  x\right)  \int_{B_{R}}\frac{V\left(  y\right)  }{\left\vert
x-y\right\vert ^{n-\beta}}\left[  \int_{B_{R}}\frac{U\left(  z\right)
\left\vert u_{k}\left(  z\right)  \right\vert ^{r}}{\left\vert y-z\right\vert
^{n-\alpha}}dz\right]  ^{1/r}dy+\eta\left(  x\right)  f_{k}\left(  x\right)  .
\end{align*}
Applying the apriori estimate to $u_{k}$, we see%
\[
\left\vert u_{k}\right\vert _{L^{q}\left(  B_{1/4}\right)  }\leq c\left(
n,p,q,r,\alpha,\beta,a,b\right)  \left(  \left\vert u_{k}\right\vert
_{L^{p}\left(  B_{1}\right)  }+\left\vert f\right\vert _{L^{q}\left(
B_{1/2}\right)  }\right)  .
\]
Now observe that%
\[
u\left(  x\right)  =T\left(  u\right)  \left(  x\right)  +\eta\left(
x\right)  f\left(  x\right)  .
\]
We see%
\begin{align*}
\left\vert u_{k}-u\right\vert _{L^{p}\left(  B_{1}\right)  }  & \leq\left\vert
T\left(  u_{k}\right)  -T\left(  u\right)  \right\vert _{L^{p}\left(
B_{1}\right)  }+\left\vert f_{k}-f\right\vert _{L^{p}\left(  B_{1}\right)  }\\
& \leq\frac{1}{2}\left\vert u_{k}-u\right\vert _{L^{p}\left(  B_{1}\right)
}+\left\vert f_{k}-f\right\vert _{L^{p}\left(  B_{1}\right)  }.
\end{align*}
Hence $\left\vert u_{k}-u\right\vert _{L^{p}\left(  B_{1}\right)  }%
\leq2\left\vert f_{k}-f\right\vert _{L^{p}\left(  B_{1}\right)  }\rightarrow0$
as $k\rightarrow\infty$. Taking a limit process in the apriori estimate for
$u_{k}$, we get the proposition.
\end{proof}

Now we are ready to derive the full regularity for the system (\ref{eq1.1}).
Such regularity under the additional assumption $p\geq1$ and $q\geq1$ was
proved in \cite{CL, L}.

\begin{proposition}
\label{propregularity}Assume $0<\alpha<n$, $\frac{\alpha}{n-\alpha}%
<p,q<\infty$, $\frac{1}{p+1}+\frac{1}{q+1}+\frac{\alpha}{n}=1$, $u\in
L_{loc}^{p+1}\left(  \mathbb{R}^{n}\right)  $ is nonnegative and does not
vanish identically. If%
\[
v\left(  x\right)  =\int_{\mathbb{R}^{n}}\frac{u\left(  y\right)  ^{p}%
}{\left\vert x-y\right\vert ^{n-\alpha}}dy,\quad u\left(  x\right)
=\int_{\mathbb{R}^{n}}\frac{v\left(  y\right)  ^{q}}{\left\vert x-y\right\vert
^{n-\alpha}}dy.
\]
Then $u,v\in C^{\infty}\left(  \mathbb{R}^{n}\right)  $. Moreover if we know
$u\in L^{p+1}\left(  \mathbb{R}^{n}\right)  $, then $u\left(  x\right)
\rightarrow0$ and $v\left(  x\right)  \rightarrow0$ as $\left\vert
x\right\vert \rightarrow\infty$.
\end{proposition}

\begin{proof}
Since $u\in L_{loc}^{p+1}\left(  \mathbb{R}^{n}\right)  $, we see $u\left(
x\right)  <\infty$ a.e. $x\in\mathbb{R}^{n}$. It follows that $v\left(
x\right)  <\infty$ a.e. $x\in\mathbb{R}^{n}$. For any $R>0$, we may find
$x_{0}\in B_{R}$ such that $v\left(  x_{0}\right)  <\infty$. This gives us
$\int_{\mathbb{R}^{n}\backslash B_{R}}\frac{u\left(  y\right)  ^{p}%
}{\left\vert x_{0}-y\right\vert ^{n-\alpha}}dy<\infty$. It follows that
$\int_{\mathbb{R}^{n}\backslash B_{R}}\frac{u\left(  y\right)  ^{p}%
}{\left\vert y\right\vert ^{n-\alpha}}dy<\infty$. Now%
\[
v\left(  x\right)  =\int_{B_{R}}\frac{u\left(  y\right)  ^{p}}{\left\vert
x-y\right\vert ^{n-\alpha}}dy+\int_{\mathbb{R}^{n}\backslash B_{R}}%
\frac{u\left(  y\right)  ^{p}}{\left\vert x-y\right\vert ^{n-\alpha}}dy,
\]
it follows from the Hardy-Littlewood-Sobolev inequality that the first term
lies in $L^{q+1}\left(  \mathbb{R}^{n}\right)  $. On the other hand, for $x\in
B_{\theta R}$ with $0<\theta<1$, we have%
\[
\int_{\mathbb{R}^{n}\backslash B_{R}}\frac{u\left(  y\right)  ^{p}}{\left\vert
x-y\right\vert ^{n-\alpha}}dy\leq\frac{1}{\left(  1-\theta\right)  ^{n-\alpha
}}\int_{\mathbb{R}^{n}\backslash B_{R}}\frac{u\left(  y\right)  ^{p}%
}{\left\vert y\right\vert ^{n-\alpha}}dy.
\]
It follows that $v\in L_{loc}^{q+1}\left(  B_{R}\right)  $. Since $R$ is
arbitrary, we have $v\in L_{loc}^{q+1}\left(  \mathbb{R}^{n}\right)  $.

Let
\begin{align*}
f_{R}\left(  x\right)   & =\int_{\mathbb{R}^{n}\backslash B_{R}}\frac{v\left(
y\right)  ^{q}}{\left\vert x-y\right\vert ^{n-\alpha}}dy,\\
g_{R}\left(  x\right)   & =\int_{\mathbb{R}^{n}\backslash B_{R}}\frac{u\left(
y\right)  ^{p}}{\left\vert x-y\right\vert ^{n-\alpha}}dy,
\end{align*}
then we know%
\begin{align*}
u\left(  x\right)   & =\int_{B_{R}}\frac{v\left(  y\right)  ^{q}}{\left\vert
x-y\right\vert ^{n-\alpha}}dy+f_{R}\left(  x\right)  ,\\
v\left(  x\right)   & =\int_{B_{R}}\frac{u\left(  y\right)  ^{p}}{\left\vert
x-y\right\vert ^{n-\alpha}}dy+g_{R}\left(  x\right)  ,
\end{align*}
and $f_{R}\in L^{p+1}\left(  B_{R}\right)  \cap L_{loc}^{\infty}\left(
B_{R}\right)  $, $g_{R}\in L^{q+1}\left(  B_{R}\right)  \cap L_{loc}^{\infty
}\left(  B_{R}\right)  $.

To continue, we observe that by symmetry, we may assume $p\geq q$, then
$p\geq\frac{n+\alpha}{n-\alpha}$ and $p-\frac{\alpha}{n}\left(  p+1\right)
\geq1$. On the other hand, it follows from $\frac{1}{p+1}+\frac{1}{q+1}%
+\frac{\alpha}{n}=1$ that $pq-1=\frac{\alpha}{n}\left(  p+1\right)  \left(
q+1\right)  $. Hence%
\[
\left[  p-\frac{\alpha}{n}\left(  p+1\right)  \right]  q-1=\frac{\alpha}%
{n}\left(  p+1\right)  >0\text{,}%
\]
and this implies $q^{-1}<p-\frac{\alpha}{n}\left(  p+1\right)  $. Choose $r$
such that%
\[
1\leq r\leq p-\frac{\alpha}{n}\left(  p+1\right)  \text{ and }q^{-1}\leq r,
\]
for example, we may take $r=p-\frac{\alpha}{n}\left(  p+1\right)  $, then%
\[
v\left(  x\right)  ^{1/r}\leq\left(  \int_{B_{R}}\frac{u\left(  y\right)
^{p}}{\left\vert x-y\right\vert ^{n-\alpha}}dy\right)  ^{1/r}+g_{R}\left(
x\right)  ^{1/r}.
\]
We have%
\begin{align*}
u\left(  x\right)   & =\int_{B_{R}}\frac{v\left(  y\right)  ^{q-r^{-1}%
}v\left(  y\right)  ^{1/r}}{\left\vert x-y\right\vert ^{n-\alpha}}%
dy+f_{R}\left(  x\right) \\
& \leq\int_{B_{R}}\frac{v\left(  y\right)  ^{q-r^{-1}}}{\left\vert
x-y\right\vert ^{n-\alpha}}\left(  \int_{B_{R}}\frac{u\left(  z\right)
^{p-r}u\left(  z\right)  ^{r}}{\left\vert y-z\right\vert ^{n-\alpha}%
}dz\right)  ^{1/r}dy+h_{R}\left(  x\right)  .
\end{align*}
Here%
\[
h_{R}\left(  x\right)  =\int_{B_{R}}\frac{v\left(  y\right)  ^{q-r^{-1}}%
g_{R}\left(  y\right)  ^{1/r}}{\left\vert x-y\right\vert ^{n-\alpha}}%
dy+f_{R}\left(  x\right)  .
\]
It follows from the fact that $g_{R}\in L^{q+1}\left(  B_{R}\right)  \cap
L_{loc}^{\infty}\left(  B_{R}\right)  $ that $h_{R}\in L^{p+1}\left(
B_{R}\right)  \cap L_{loc}^{\overline{q}}\left(  B_{R}\right)  $ for all
$\overline{q}<\infty$. Let%
\[
a=\frac{p+1}{p-r},\quad b=\frac{q+1}{q-r^{-1}},
\]
then calculation shows $\frac{1}{ra}+\frac{1}{b}=\frac{\alpha}{rn}%
+\frac{\alpha}{n}$, moreover we have%
\[
\frac{r}{p+1}+\frac{1}{a}=\frac{p}{p+1}<1,
\]
and%
\[
\frac{1}{a}-\frac{\alpha}{n}=\frac{p-\frac{\alpha}{n}\left(  p+1\right)
-r}{p+1}\geq0.
\]
Hence for any $p+1<\overline{q}<\infty$, when $R$ is small enough, it follows
from Proposition \ref{proplocal} (by choosing $\alpha$, $\beta$, $p $, $q$,
$r$, $a$, $b$, $u$, $U$, $V$ and $f$ in Proposition \ref{proplocal} as
$\alpha$, $\alpha$, $p+1$, $\overline{q}$, $r$, $\frac{p+1}{p-r}$, $\frac
{q+1}{q-r^{-1}}$, $u$, $u^{p-r}$, $v^{q-r^{-1}}$ and $h_{R}$ respectively)
that $u\in L^{\overline{q}}\left(  B_{R/4}\right)  $. Since every point may be
viewed as a center, we see $u\in L_{loc}^{\overline{q}}\left(  \mathbb{R}%
^{n}\right)  $. This implies $v\in L_{loc}^{\infty}\left(  \mathbb{R}%
^{n}\right)  $ and then $u\in L_{loc}^{\infty}\left(  \mathbb{R}^{n}\right)
$. Now observe that $f_{R},g_{R}\in C^{\infty}\left(  B_{R}\right)  $, it
follows from the usual bootstrap method that $u,v\in C^{\infty}\left(
\mathbb{R}^{n}\right)  $. The fact $u,v\in L^{\infty}\left(  \mathbb{R}%
^{n}\right)  $ under the assumption $u\in L^{p+1}\left(  \mathbb{R}%
^{n}\right)  $ follows from carefully going through the above argument and
applying Holder's inequality when needed. Note that%
\[
u=\frac{\chi_{B_{1}}\left(  x\right)  }{\left\vert x\right\vert ^{n-\alpha}%
}\ast v^{q}+\frac{\chi_{\mathbb{R}^{n}\backslash B_{1}}\left(  x\right)
}{\left\vert x\right\vert ^{n-\alpha}}\ast v^{q}.
\]
By interpolation we know $v\in L^{s}\left(  \mathbb{R}^{n}\right)  $ for all
$q+1\leq s\leq\infty$, hence $v^{q}\in L^{s}$ for $\frac{q+1}{q}\leq
s\leq\infty$. Since $\frac{q+1}{q}<\frac{n}{\alpha}$, it follows form the fact
$\frac{\chi_{B_{1}}\left(  x\right)  }{\left\vert x\right\vert ^{n-\alpha}}\in
L^{\frac{n}{n-\alpha}-\varepsilon}\left(  \mathbb{R}^{n}\right)  $ and
$\frac{\chi_{\mathbb{R}^{n}\backslash B_{1}}\left(  x\right)  }{\left\vert
x\right\vert ^{n-\alpha}}\in L^{\frac{n}{n-\alpha}+\varepsilon}\left(
\mathbb{R}^{n}\right)  $ for $\varepsilon>0$ small that $u\left(  x\right)
\rightarrow0$ as $\left\vert x\right\vert \rightarrow\infty$. The fact
$v\left(  x\right)  \rightarrow0$ as $\left\vert x\right\vert \rightarrow
\infty$ follows similarly.
\end{proof}

\section{All solutions are radial\label{sec3}}

In this section, we will use the integral form of the method of moving plane
(developed in \cite{CLO2}) to prove the radial symmetry of solutions to the
integral system. Such radial property was derived in \cite{CLO1,CLO2,L} under
the further assumptions that both $p$ and $q$ are at least $1$. Our approach
works for both this case and the case when $p$ or $q$ is strictly less than
$1$. We will need the following basic inequality: assume $0<\theta\leq1$,
$a\geq b\geq0$, $c\geq0$, then%
\[
\left(  a+c\right)  ^{\theta}-\left(  b+c\right)  ^{\theta}\leq a^{\theta
}-b^{\theta}.
\]
Indeed, for $x\geq0$, let $f\left(  x\right)  =\left(  a+x\right)  ^{\theta
}-\left(  b+x\right)  ^{\theta}$, then for $x>0$, $f^{\prime}\left(  x\right)
=\theta\left(  a+x\right)  ^{\theta-1}-\theta\left(  b+x\right)  ^{\theta
-1}\leq0$. The inequality follows.

For $\xi\in\mathbb{R}^{m}$ and $s>0$, we denote%
\[
\left\vert \xi\right\vert _{l^{s}}=\left(  \sum_{i=1}^{m}\left\vert \xi
_{i}\right\vert ^{s}\right)  ^{1/s}.
\]

\begin{proof}
[Proof of Theorem \ref{thmmain}]By Proposition \ref{propregularity}, we know
$u,v\in C^{\infty}\left(  \mathbb{R}^{n}\right)  $, $u\left(  x\right)
\rightarrow0$ and $v\left(  x\right)  \rightarrow0$ as $\left\vert
x\right\vert \rightarrow\infty$. It follows from Hardy-Littlewood-Sobolev
inequality that $v\in L^{q+1}\left(  \mathbb{R}^{n}\right)  $. Without losing
of generality, we may assume $p\geq q $, then we know $p\geq\frac{n+\alpha
}{n-\alpha}$ and $p>q^{-1}$. Hence we may find a number $r$ such that $1\leq
r<p$ and $q^{-1}\leq r$.

For $\lambda\in\mathbb{R}$, we denote $H_{\lambda}=\left\{  x\in\mathbb{R}%
^{n}:x_{1}<\lambda\right\}  $. For $x=\left(  x_{1},x^{\prime}\right)
\in\mathbb{R}^{n}$, let $x_{\lambda}=\left(  2\lambda-x_{1},x^{\prime}\right)
$. We also denote $u_{\lambda}\left(  x\right)  =u\left(  x_{\lambda}\right)
$, $v_{\lambda}\left(  x\right)  =v\left(  x_{\lambda}\right)  $,%
\begin{align*}
\mathcal{B}_{\lambda}^{u}  & =\left\{  x\in H_{\lambda}:u_{\lambda}\left(
x\right)  >u\left(  x\right)  \right\}  ,\\
\mathcal{B}_{\lambda}^{v}  & =\left\{  x\in H_{\lambda}:v_{\lambda}\left(
x\right)  >v\left(  x\right)  \right\}  .
\end{align*}
Note that by a change of variable, we have%
\begin{align*}
u\left(  x\right)   & =\int_{\mathbb{R}^{n}}\frac{v\left(  y\right)  ^{q}%
}{\left\vert x-y\right\vert ^{n-\alpha}}dy\\
& =\int_{H_{\lambda}}\frac{v\left(  y\right)  ^{q}}{\left\vert x-y\right\vert
^{n-\alpha}}dy+\int_{H_{\lambda}}\frac{v\left(  y_{\lambda}\right)  ^{q}%
}{\left\vert x_{\lambda}-y\right\vert ^{n-\alpha}}dy.
\end{align*}
Hence%
\[
u\left(  x_{\lambda}\right)  =\int_{H_{\lambda}}\frac{v\left(  y_{\lambda
}\right)  ^{q}}{\left\vert x-y\right\vert ^{n-\alpha}}dy+\int_{H_{\lambda}%
}\frac{v\left(  y\right)  ^{q}}{\left\vert x_{\lambda}-y\right\vert
^{n-\alpha}}dy.
\]
This implies%
\begin{align*}
& u\left(  x_{\lambda}\right)  -u\left(  x\right) \\
& =\int_{H_{\lambda}}\left(  v\left(  y_{\lambda}\right)  ^{q}-v\left(
y\right)  ^{q}\right)  \left(  \frac{1}{\left\vert x-y\right\vert ^{n-\alpha}%
}-\frac{1}{\left\vert x_{\lambda}-y\right\vert ^{n-\alpha}}\right)  dy.
\end{align*}
In particular, for $x\in\mathcal{B}_{\lambda}^{u}$, we have%
\begin{align*}
0  & \leq u\left(  x_{\lambda}\right)  -u\left(  x\right) \\
& \leq\int_{\mathcal{B}_{\lambda}^{v}}\left(  v\left(  y_{\lambda}\right)
^{q}-v\left(  y\right)  ^{q}\right)  \left(  \frac{1}{\left\vert
x-y\right\vert ^{n-\alpha}}-\frac{1}{\left\vert x_{\lambda}-y\right\vert
^{n-\alpha}}\right)  dy\\
& \leq\int_{\mathcal{B}_{\lambda}^{v}}\left(  \left(  v\left(  y_{\lambda
}\right)  ^{1/r}\right)  ^{qr}-\left(  v\left(  y\right)  ^{1/r}\right)
^{qr}\right)  \frac{1}{\left\vert x-y\right\vert ^{n-\alpha}}dy\\
& \leq qr\int_{\mathcal{B}_{\lambda}^{v}}v\left(  y_{\lambda}\right)
^{q-r^{-1}}\left(  v\left(  y_{\lambda}\right)  ^{1/r}-v\left(  y\right)
^{1/r}\right)  \frac{1}{\left\vert x-y\right\vert ^{n-\alpha}}dy.
\end{align*}
It follows from Hardy-Littlewood-Sobolev inequality that%
\begin{align*}
& \left\vert u_{\lambda}-u\right\vert _{L^{p+1}\left(  \mathcal{B}_{\lambda
}^{u}\right)  }\\
& \leq c\left(  n,\alpha,q,r\right)  \left\vert v_{\lambda}^{q-r^{-1}}\left(
v_{\lambda}^{1/r}-v^{1/r}\right)  \right\vert _{L^{\frac{q+1}{q}}\left(
\mathcal{B}_{\lambda}^{v}\right)  }\\
& \leq c\left(  n,\alpha,q,r\right)  \left\vert v_{\lambda}^{q-r^{-1}%
}\right\vert _{L^{\frac{q+1}{q-r^{-1}}}\left(  \mathcal{B}_{\lambda}%
^{v}\right)  }\left\vert v_{\lambda}^{1/r}-v^{1/r}\right\vert _{L^{\left(
q+1\right)  r}\left(  \mathcal{B}_{\lambda}^{v}\right)  }\\
& =c\left(  n,\alpha,q,r\right)  \left\vert v_{\lambda}\right\vert
_{L^{q+1}\left(  \mathcal{B}_{\lambda}^{v}\right)  }^{q-r^{-1}}\left\vert
v_{\lambda}^{1/r}-v^{1/r}\right\vert _{L^{\left(  q+1\right)  r}\left(
\mathcal{B}_{\lambda}^{v}\right)  }%
\end{align*}
On the other hand, for $x\in\mathcal{B}_{\lambda}^{v}$, we have%
\begin{align*}
v\left(  x_{\lambda}\right)   & =\int_{\mathcal{B}_{\lambda}^{u}}%
\frac{u\left(  y_{\lambda}\right)  ^{p}}{\left\vert x-y\right\vert ^{n-\alpha
}}dy+\int_{\mathcal{B}_{\lambda}^{u}}\frac{u\left(  y\right)  ^{p}}{\left\vert
x_{\lambda}-y\right\vert ^{n-\alpha}}dy\\
& +\int_{H_{\lambda}\backslash\mathcal{B}_{\lambda}^{u}}\frac{u\left(
y_{\lambda}\right)  ^{p}}{\left\vert x-y\right\vert ^{n-\alpha}}%
dy+\int_{H_{\lambda}\backslash\mathcal{B}_{\lambda}^{u}}\frac{u\left(
y\right)  ^{p}}{\left\vert x_{\lambda}-y\right\vert ^{n-\alpha}}dy\\
& \leq\int_{\mathcal{B}_{\lambda}^{u}}\frac{u\left(  y_{\lambda}\right)  ^{p}%
}{\left\vert x-y\right\vert ^{n-\alpha}}dy+\int_{\mathcal{B}_{\lambda}^{u}%
}\frac{u\left(  y\right)  ^{p}}{\left\vert x_{\lambda}-y\right\vert
^{n-\alpha}}dy\\
& +\int_{H_{\lambda}\backslash\mathcal{B}_{\lambda}^{u}}\frac{u\left(
y\right)  ^{p}}{\left\vert x-y\right\vert ^{n-\alpha}}dy+\int_{H_{\lambda
}\backslash\mathcal{B}_{\lambda}^{u}}\frac{u\left(  y_{\lambda}\right)  ^{p}%
}{\left\vert x_{\lambda}-y\right\vert ^{n-\alpha}}dy.
\end{align*}
Since
\begin{align*}
v\left(  x\right)   & =\int_{\mathcal{B}_{\lambda}^{u}}\frac{u\left(
y\right)  ^{p}}{\left\vert x-y\right\vert ^{n-\alpha}}dy+\int_{\mathcal{B}%
_{\lambda}^{u}}\frac{u\left(  y_{\lambda}\right)  ^{p}}{\left\vert x_{\lambda
}-y\right\vert ^{n-\alpha}}dy\\
& +\int_{H_{\lambda}\backslash\mathcal{B}_{\lambda}^{u}}\frac{u\left(
y\right)  ^{p}}{\left\vert x-y\right\vert ^{n-\alpha}}dy+\int_{H_{\lambda
}\backslash\mathcal{B}_{\lambda}^{u}}\frac{u\left(  y_{\lambda}\right)  ^{p}%
}{\left\vert x_{\lambda}-y\right\vert ^{n-\alpha}}dy,
\end{align*}
it follows that%
\begin{align*}
0  & \leq v\left(  x_{\lambda}\right)  ^{1/r}-v\left(  x\right)  ^{1/r}\\
& \leq\left(  \int_{\mathcal{B}_{\lambda}^{u}}\frac{u\left(  y_{\lambda
}\right)  ^{p}}{\left\vert x-y\right\vert ^{n-\alpha}}dy+\int_{\mathcal{B}%
_{\lambda}^{u}}\frac{u\left(  y\right)  ^{p}}{\left\vert x_{\lambda
}-y\right\vert ^{n-\alpha}}dy\right)  ^{1/r}\\
& -\left(  \int_{\mathcal{B}_{\lambda}^{u}}\frac{u\left(  y\right)  ^{p}%
}{\left\vert x-y\right\vert ^{n-\alpha}}dy+\int_{\mathcal{B}_{\lambda}^{u}%
}\frac{u\left(  y_{\lambda}\right)  ^{p}}{\left\vert x_{\lambda}-y\right\vert
^{n-\alpha}}dy\right)  ^{1/r}\\
& =\left(  \int_{\mathcal{B}_{\lambda}^{u}}\left\vert \left(  \frac{u\left(
y_{\lambda}\right)  ^{p/r}}{\left\vert x-y\right\vert ^{\left(  n-\alpha
\right)  /r}},\frac{u\left(  y\right)  ^{p/r}}{\left\vert x_{\lambda
}-y\right\vert ^{\left(  n-\alpha\right)  /r}}\right)  \right\vert _{l^{r}%
}^{r}dy\right)  ^{1/r}\\
& -\left(  \int_{\mathcal{B}_{\lambda}^{u}}\left\vert \left(  \frac{u\left(
y\right)  ^{p/r}}{\left\vert x-y\right\vert ^{\left(  n-\alpha\right)  /r}%
},\frac{u\left(  y_{\lambda}\right)  ^{p/r}}{\left\vert x_{\lambda
}-y\right\vert ^{\left(  n-\alpha\right)  /r}}\right)  \right\vert _{l^{r}%
}^{r}dy\right)  ^{1/r}\\
& \leq\left(  \int_{\mathcal{B}_{\lambda}^{u}}\left\vert \left(
\frac{u\left(  y_{\lambda}\right)  ^{p/r}-u\left(  y\right)  ^{p/r}%
}{\left\vert x-y\right\vert ^{\left(  n-\alpha\right)  /r}},\frac{u\left(
y\right)  ^{p/r}-u\left(  y_{\lambda}\right)  ^{p/r}}{\left\vert x_{\lambda
}-y\right\vert ^{\left(  n-\alpha\right)  /r}}\right)  \right\vert _{l^{r}%
}^{r}dy\right)  ^{1/r}\\
& \leq2\left(  \int_{\mathcal{B}_{\lambda}^{u}}\frac{\left(  u_{\lambda
}\left(  y\right)  ^{p/r}-u\left(  y\right)  ^{p/r}\right)  ^{r}}{\left\vert
x-y\right\vert ^{n-\alpha}}dy\right)  ^{1/r}\\
& \leq\frac{2p}{r}\left(  \int_{\mathcal{B}_{\lambda}^{u}}\frac{u_{\lambda
}\left(  y\right)  ^{p-r}\left(  u_{\lambda}\left(  y\right)  -u\left(
y\right)  \right)  ^{r}}{\left\vert x-y\right\vert ^{n-\alpha}}dy\right)
^{1/r}.
\end{align*}
It follows from Hardy-Littlewood-Sobolev inequality that%
\begin{align*}
& \left\vert v_{\lambda}^{1/r}-v^{1/r}\right\vert _{L^{\left(  q+1\right)
r}\left(  \mathcal{B}_{\lambda}^{v}\right)  }\\
& \leq\frac{2p}{r}\left\vert \int_{\mathcal{B}_{\lambda}^{u}}\frac{u_{\lambda
}\left(  y\right)  ^{p-r}\left(  u_{\lambda}\left(  y\right)  -u\left(
y\right)  \right)  ^{r}}{\left\vert x-y\right\vert ^{n-\alpha}}dy\right\vert
_{L^{q+1}\left(  \mathcal{B}_{\lambda}^{v}\right)  }^{1/r}\\
& \leq c\left(  n,\alpha,p,r\right)  \left\vert u_{\lambda}^{p-r}\left(
u_{\lambda}-u\right)  ^{r}\right\vert _{L^{\frac{p+1}{p}}\left(
\mathcal{B}_{\lambda}^{u}\right)  }^{1/r}\\
& \leq c\left(  n,\alpha,p,r\right)  \left\vert u_{\lambda}^{p-r}\right\vert
_{L^{\frac{p+1}{p-r}}\left(  \mathcal{B}_{\lambda}^{u}\right)  }%
^{1/r}\left\vert \left(  u_{\lambda}-u\right)  ^{r}\right\vert _{L^{\left(
p+1\right)  /r}\left(  \mathcal{B}_{\lambda}^{u}\right)  }^{1/r}\\
& =c\left(  n,\alpha,p,r\right)  \left\vert u_{\lambda}\right\vert
_{L^{p+1}\left(  \mathcal{B}_{\lambda}^{u}\right)  }^{\frac{p-r}{r}}\left\vert
u_{\lambda}-u\right\vert _{L^{p+1}\left(  \mathcal{B}_{\lambda}^{u}\right)  }.
\end{align*}
Hence we have%
\begin{align*}
& \left\vert u_{\lambda}-u\right\vert _{L^{p+1}\left(  \mathcal{B}_{\lambda
}^{u}\right)  }\\
& \leq c\left(  n,\alpha,p,q,r\right)  \left\vert u_{\lambda}\right\vert
_{L^{p+1}\left(  \mathcal{B}_{\lambda}^{u}\right)  }^{\frac{p-r}{r}}\left\vert
v_{\lambda}\right\vert _{L^{q+1}\left(  \mathcal{B}_{\lambda}^{v}\right)
}^{q-r^{-1}}\left\vert u_{\lambda}-u\right\vert _{L^{p+1}\left(
\mathcal{B}_{\lambda}^{u}\right)  }\\
& =c\left(  n,\alpha,p,q,r\right)  \left\vert u\right\vert _{L^{p+1}\left(
2\lambda e_{1}-\mathcal{B}_{\lambda}^{u}\right)  }^{\frac{p-r}{r}}\left\vert
v\right\vert _{L^{q+1}\left(  2\lambda e_{1}-\mathcal{B}_{\lambda}^{v}\right)
}^{q-r^{-1}}\left\vert u_{\lambda}-u\right\vert _{L^{p+1}\left(
\mathcal{B}_{\lambda}^{u}\right)  }\\
& \leq c\left(  n,\alpha,p,q,r\right)  \left\vert u\right\vert _{L^{p+1}%
\left(  2\lambda e_{1}-\mathcal{B}_{\lambda}^{u}\right)  }^{\frac{p-r}{r}%
}\left\vert v\right\vert _{L^{q+1}\left(  \mathbb{R}^{n}\right)  }^{q-r^{-1}%
}\left\vert u_{\lambda}-u\right\vert _{L^{p+1}\left(  \mathcal{B}_{\lambda
}^{u}\right)  }.
\end{align*}
Here $e_{1}=\left(  1,0,\cdots,0\right)  $.

After these preparations, we will use the method of moving planes to prove the
radial symmetry of the solutions.

First, we have to show it is possible to start. Indeed, for $\lambda$ large
enough, we know $\left\vert u\right\vert _{L^{p+1}\left(  2\lambda
e_{1}-\mathcal{B}_{\lambda}^{u}\right)  }$ can be arbitrary small, this
implies that%
\[
\left\vert u_{\lambda}-u\right\vert _{L^{p+1}\left(  \mathcal{B}_{\lambda}%
^{u}\right)  }\leq\frac{1}{2}\left\vert u_{\lambda}-u\right\vert
_{L^{p+1}\left(  \mathcal{B}_{\lambda}^{u}\right)  },
\]
and hence $\left\vert u_{\lambda}-u\right\vert _{L^{p+1}\left(  \mathcal{B}%
_{\lambda}^{u}\right)  }=0$. It follows that $\mathcal{B}_{\lambda}%
^{u}=\emptyset$ when $\lambda$ is large enough.

Next we let $\lambda_{0}=\inf\left\{  \lambda\in\mathbb{R}:\mathcal{B}%
_{\lambda^{\prime}}^{u}=\emptyset\text{ for all }\lambda^{\prime}\geq
\lambda\right\}  $. It follows from the fact $u\left(  x\right)  \rightarrow0$
as $\left\vert x\right\vert \rightarrow\infty$ and $u\left(  x\right)  >0$ for
all $x$ that $\lambda_{0}$ must be a finite number. It follows from the
definition of $\lambda_{0}$ that $u_{\lambda_{0}}\left(  x\right)  \leq
u\left(  x\right)  $ for $x\in H_{\lambda_{0}}$. We claim that $u_{\lambda
_{0}}=u$. If this is not the case, then since%
\[
v_{\lambda_{0}}\left(  x\right)  -v\left(  x\right)  =\int_{H_{\lambda_{0}}%
}\left(  u_{\lambda_{0}}\left(  y\right)  ^{p}-u\left(  y\right)  ^{p}\right)
\left(  \frac{1}{\left\vert x-y\right\vert ^{n-\alpha}}-\frac{1}{\left\vert
x_{\lambda_{0}}-y\right\vert ^{n-\alpha}}\right)  dy
\]
and%
\[
u_{\lambda_{0}}\left(  x\right)  -u\left(  x\right)  =\int_{H_{\lambda_{0}}%
}\left(  v_{\lambda_{0}}\left(  y\right)  ^{q}-v\left(  y\right)  ^{q}\right)
\left(  \frac{1}{\left\vert x-y\right\vert ^{n-\alpha}}-\frac{1}{\left\vert
x_{\lambda_{0}}-y\right\vert ^{n-\alpha}}\right)  dy,
\]
we see $u_{\lambda_{0}}\left(  x\right)  <u\left(  x\right)  $ for $x\in
H_{\lambda_{0}}$. This implies $\chi_{2\lambda e_{1}-\mathcal{B}_{\lambda}%
^{u}}\rightarrow0$ $a.e.$ as $\lambda\uparrow\lambda_{0}$. It follows that
$\left\vert u\right\vert _{L^{p+1}\left(  2\lambda e_{1}-\mathcal{B}_{\lambda
}^{u}\right)  }\rightarrow0$ as $\lambda\uparrow\lambda_{0}$. Hence%
\[
\left\vert u_{\lambda}-u\right\vert _{L^{p+1}\left(  \mathcal{B}_{\lambda}%
^{u}\right)  }\leq\frac{1}{2}\left\vert u_{\lambda}-u\right\vert
_{L^{p+1}\left(  \mathcal{B}_{\lambda}^{u}\right)  }%
\]
when $\lambda$ is close to $\lambda_{0}$. This implies $\mathcal{B}_{\lambda
}^{u}=\emptyset$ for $\lambda$ close to $\lambda_{0}$ and it contradicts with
the choice of $\lambda_{0}$. Hence when the moving plane process stops, we
must have symmetry. Moreover, $u_{\lambda}\left(  x\right)  <u\left(
x\right)  $ for $x\in H_{\lambda}$ when $\lambda>\lambda_{0}$. Indeed, for any
$\lambda>\lambda_{0}$, we can not have $u_{\lambda}=u$ because otherwise $u$
is periodic in the first direction and can not lie in $L^{p+1}$. Hence
$u_{\lambda}<u$ in $H_{\lambda}$.

By translation, we may assume $u\left(  0\right)  =\max_{x\in\mathbb{R}^{n}%
}u\left(  x\right)  $, then it follows that the moving plane process from any
direction must stop at the origin, hence $u$ must be radial symmetric and
strictly decreasing in the radial direction. It follows from the equation that
$v$ has the same properties.
\end{proof}

\end{document}